\documentstyle[12pt,amssymb]{article}
\begin{document}

\newtheorem{thm}{Theorem}
\newtheorem{lmm}{Lemma}
\newtheorem{dfn}{Definition}
\newtheorem{rmk}{Remark}
\newtheorem{prp}{Proposition}
\newtheorem{exa}{Example}
\newtheorem{cor}{Corollary}

\title{MEROMORPHIC BRAIDED CATEGORY ARISING IN QUANTUM AFFINE ALGEBRAS}

\author{Y. Soibelman }

\maketitle

\section{Introduction}
Main result of this article consists in construction
 of what is called meromorphic
braided (or tensor) category, see [So]. 
It  arises in the representation
theory of quantum affine algebras and lives on an elliptic curve. 

There are many papers devoted to different aspects
 of finite-dimensional
representations of quantum affine algebras.
 Surprisinly few of them consider
categorical picture. Most of the authors
 study  
 irreducible representations. 
Some fundamental questions remain unanswered.
For example the fact that 
the universal $R$-matrix is  meromorphic for $any$
two finite-dimensional
 representations was not proved (to my knowledge)
before [KS]. This fact is crucial 
for construction of meromorphic braided structure
on the category  of finite-dimensional
 representations. The latter contains an
interesting subcategory with objects 
naturally ``localized'' on an elliptic curve.

 In this article we consider the simplest example related
to the quantum affine algebra $sl(2)$.
 Our constructions remind 
``chiral'' objects from [BD]: we  derive the 
braiding considering infinitesimal neighbourhood
of the diagonal in the square of an elliptic curve.

The paper is organized as follows. 
Next two sections contain recollections from [So].
Last section contains construction and 
discussion. It can be read  independently of
the other part of the paper.
In the last section we mainly
 concentrate on the case of quantum affine $sl(2)$.
It is explained at the very end of the paper 
how main construction can be generalized
to higher ranks if one uses results of [FR].

$Acknowledgements$. Constructions of Section 4
 have arisen from a number of
 discussions with Joseph Bernstein. 
I am also greatful to Ed Frenkel 
for the discussions about $q$-characters 
and simple modules over affine quantum
algebras.

\section{Pseudo-braided categories}

In this section we are going to recall (in a slightly revised form)
main definitions from [So].

\subsection{}

We denote by ${\cal T}$ the class of all planar trees, and by ${\cal T}(n)$
the subclass of trees having $n$  tails (see [GK], [KM] about terminology).
Then ${\cal T}(n)$ is a category (morphisms are identities and
contraction of  edges). For a given tree we  orient
edges and tails in such a way that there is a unique vertex with
only one outgoing tail (this vertex is called the root of tree).
We fix such an orientation for each tree. We also fix a numeration
of tails : for each $T \in {\cal T}(n)$ they are numbered from $1 $ to $n$
in such a way that the only outgoing tail is numbered by $n$.

An additional structure on ${\cal T}$ is given by the gluing operation:
if $T \in {\cal T}(n+1), T_i \in {\cal T}(k_i+1), i=1,...,n+1$ then one can
construct a new tree $T(T_1,...,T_n) \in {\cal T}(k_1+...+k_n+1)$ by gluing
outgoing tail of $T_i$ to the $i$th tail of $T$. The orientation of edges
and numeration of tails for the new tree are defined in the natural way.
The gluing operation is associative, and hence ${\cal T}$ becomes a strict
monoidal 2-operad . The role of a unit object is played by the only tree
$e \in {\cal T}(1)$.

Let  ${\cal S}=(S_T,{\cal O}_{S_T})$ be a family of ringed spaces
 parametrized by trees 
from ${\cal T}$. We say that ${\cal S}$ is a monoidal
operad of spaces if the following conditions are satisfied:

1) for any morphism $f: T^{\prime} \rightarrow T \in {\cal T}(n)$
we are given a morphism of ringed spaces
 $l_f: S_{T^{\prime}} \rightarrow S_T$.

2) For a gluing operation of trees 
$T \times T_1 \times...\times T_n \rightarrow
T(T_1,...,T_n)$ we are given a morphism of ringed spaces
(operadic composition) $\gamma:
S_T \times S_{T_1} \times ... \times S_{T_n} \rightarrow S_{T(T_1,...,T_n)}$
which is strictly associative with respect to the gluing of trees.
It is also assumed to be functorial with respect to the morphisms
of trees.

In particular if we put $A=\Gamma(S_e,{\cal O}_{S_e})$,
then all $\Gamma(S_T,{\cal O}_{S_T})$ become 
$A-A^{\otimes n}$ bimodules, where $T \in {\cal T}(n)$.

\subsection{}
Let $X$ be a set, ${\cal A}$  a class. Its elements are called objects.
A family of objects of ${\cal A}$ parametrized by $X$ (or simply
$X$-family) is an element of $\prod_X {\cal A}= {\cal A}^X$.
Suppose that ${\cal A}$ is  a category. We keep the same notation
for the class of objects of ${\cal A}$.
Then $X$-families form a category ${\cal F}_X$ with 
$Hom_{{\cal F}_X}(M,N)=
\prod_X Hom_{\cal A}(M_x,N_x)$.
If $f: X \rightarrow Y$ is a map the there is a pull-back
functor $f^{*}: {\cal F}_Y \rightarrow {\cal F}_X$.

If $(X,{\cal O}_X)$ is a ringed space then we understand
an $X$-family as a family of objects  which
are ${\cal O}_{X,x}$-modules, where $O_{X,x}$ denotes
the fiber at $x \in X$ . The pull-back
functors exist  in the case of ringed spaces as well.
We do not impose any continuity condition on the family
at the moment.

Suppose that  ${\cal S}$ is a monoidal operad of spaces
as in the  previous subsection. We suppose that we are
given a field $k$ , and that all sheaves of rings are
in fact $k$-algebras and all morphisms of sheaves of
rings are morphisms of $k$-algebras.

\begin{dfn}

A pseudo-monoidal category ${\cal C}$
of ringed spaces is given by the following data:

a) a class ${\cal C}$ called class of objects;

b) for  every $T \in {\cal T}(n )$ , a sequence
${\{X_i \}}, 1 \le i \le n$ of objects of ${\cal C}$
and an object $Y \in {\cal C}$ a family of $k$-vector
spaces $P_T(\{X_i\},Y)$ over $S_T$ (operations from
$\{X_i \}$ to $Y$);

c) for every morphism $f: T^{\prime} \rightarrow T$
a morphism $\phi_f: P_{T^{\prime}} (\{X_i \} ,Y) \rightarrow 
(l_f)^{*}P_T (\{X_i \},Y)$ of families over $S_{T^{\prime}}$;

d) for an operadic composition of spaces

 $\gamma: S_T \times \prod_i S_{T_i} \rightarrow S_{T(T_1,...,T_n)}$
we are given a morphism of families on $S_T \times \prod_i S_{T_i}$
(composition of operations):

$\Phi_{\gamma}: P_T(\{X_i \},Y) \times \prod_i P_{T_i} (\{ K_{j_i} \}, X_i)
\rightarrow \gamma^{*} P_{T(T_1,...,T_n)} (\{K_j \}, Y)$ , 

$\Phi_{\gamma}(\phi , (\psi_i))= \phi (\psi_i)$.

 Here $\{X_i \}$ are parametrized by the tails
 of $T$, $\{ K_j \}$ are parametrized
by the tails of $T(T_1,...,T_n)$ and $\{ K_{j_i} \}$ 
corresponds to the subsequence
of objects parametrized by the tails of $T_i$;

e) for every object $X \in {\cal C}$ there exists 
a family ${\bf 1}_X \in P_e( \{ X \},X)$
on $S_e$ such that $\phi ({\bf 1}_{X_i})= \gamma^{*}(\phi)$
and ${\bf 1}_X (\phi)= \eta^{*}(\phi)$. 

Here $\phi \in P_T(\{ X_i \},X)$ and we use
 the notations $\gamma$ and $\eta$
for the natural operadic morphisms of spaces (see d)).

The composition morphisms $\Phi_{\gamma}$ from d)
 are required to be associative
in the following sense: 

$\Phi_{\gamma}( id \times \Phi_{\delta})=\Phi_{\delta}(\Phi_{\gamma}
\times id)$ 

(when applied to the corresponding families).

\end{dfn}

Sometimes we will call ${\cal C}$
 a pseudo-monoidal category over ${\cal S}$.

Suppose our spaces are $k$-schemes.
 We call a pseudo-monoidal category
$algebraic$ if all families $P_T(\{X_i \}, Y)$ 
are quasi-coherent sheaves
and the corresponding morphisms
 are morphisms of quasi-coherent sheaves.

Suppose that  $k={\bf C}$ and our spaces
 are complex analytic. We call
a pseudo-monoidal category $analytic$
 if the families $P_T(\{X_i \},Y)$
are complex analytic sheaves
 and the corresponding morphisms
are morphisms of complex analytic sheaves.

Similarly, if these complex analytic spaces
 are irreducible, we can
speak about $meromorphic$ pseudo-monoidal 
categories (in this
case $P_T$ are analytic families on dense 
subsets and morphisms
can be extended meromorphically to $S_T$).

In the case of schemes we obtain $rational$ 
pseudo-monoidal
categories.

We are going to use this terminology without further discussion.

\subsection{}

Here we recall what is it a representable
 pseudo-monoidal structure.
Let $S$ be a topological space, ${\cal C} $ 
 a category, $n \ge 1$ an integer.

We denote by $Funct_S({\cal C},n)$ the
 sheaf of categories on $S$ such
that for an open $U$ in $S$ we have
 $Funct_S({\cal C},n)(U)$= category
of families of functors
 $\{F_x \}_{x \in U}, F_x: {\cal C}^n \rightarrow
{\cal C}$. We will denote this sheaf
 of categories by $Funct_S(n)$ if 
it will not lead to a confusion.

Suppose that in addition we have 
a sequence of topological spaces
$\{ S_i \}, 1\le i\le n$ , and a sequence
 of positive integers
$k_1,...,k_n$ . Then we have a morphism
of sheaves of categories $Funct_S(n) \times \prod_iFunct_{S_i}
\rightarrow Funct_{S \times \prod_iS_i}(k_1+...+k_n)$ (sheaves on 
$S \times \prod_iS_i$) such that $\{F_s\}\times\prod_iF_{s_i}^i
\rightarrow F_s(F^1_{s_1},...,F^n_{s_n})$.

If ${\cal C}$ is a $k$-linear category then for a family $\{F_x \}
\in Funct_S(n)(S)$ and a sequence $\{X_i \}, 1\le i\le n$ of objects
of ${\cal C}$ and an object $Y \in {\cal C}$ we have a family of
vectors spaces $Hom_{{\cal C}}(F_x(\{X_i \}), Y)$ on $S$.

Suppose now that we have a monoidal operad of spaces ${\cal S}=
(S_T)$ as before. Suppose that for every $T\in {\cal T}(n)$ we 
are given a family of functors $\{F_x^T\}_{x \in S_T} \in
Funct_{S_T}(n)$. We remark that there is a natural operadic
composition $Funct_{S_T}(n)\times \prod_i Funct_{S_{T_i}}(k_i)
\rightarrow Funct_{S_T \times \prod_iS_{T_i}}(k_1+...+k_n) \rightarrow
\gamma^{*}(Funct_{S_{T(T_1,...,T_n)}}(k_1+...+k_n)$ of 
sheaves of categories.

The last arrow corresponds to the operadic
 composition $\gamma$ on 
the spaces which induces the obvious pull-back on the families.
Suppose that we are given a pseudo-monoidal 
category on ${\cal S}$ as
before.
\begin{dfn}
We say that it is representable if for every
 $T \in {\cal T}(n)$ there
exists a family $F_x^T \in Funct_{S_T}(n)(S_T)$
 and an isomorphism
of families of $k$- vector spaces
 $P_T(\{X_i \},Y) \rightarrow
Hom_{{\cal C}}(F_x^T(\{X_i\}),Y)$ which is compatible with
the operadic composition on both families.

It is also required that for the only tree $e \in {\cal T}(1)$ we
have: $Hom_{{\cal C}}(F_x^e(Y),Y)$ corresponds to ${\bf 1}_Y$ 
under this isomorphism.

\end{dfn}

If ${\cal C}$ is an algebraic or analytic pseudo-monoidal category
then representability is understood in the corresponding category.
In particular it is assumed that families of functors define families
of $Hom$'s in the corresponding category, and the isomorphisms
of families must be isomorphisms of quasi-coherent sheaves or
analytic sheaves. Similarly one defines representable rational
and meromorphic pseudo-monoidal categories. We skip "pseudo" in the
case when all morphisms in the Definition 1c) are isomorphisms.
In this way we obtain for example the notion of meromorphic monoidal
category discussed in [So].

It is easy to define the notion of a functor between two pseudo-monoidal
categories ${\cal A}$ and ${\cal B}$ which live over different
operads of spaces,say, ${\cal S}=(S_T)$ and ${\cal R}=(R_T)$ respectively.
It consists of  morphisms of ringed spaces $h_T:S_T \rightarrow
R_T$ which produce a morphism of the monoidal operads of spaces,
 of the mapping of objects  $F: {\cal A} \rightarrow {\cal B}$, and
of the morphisms 
of families $l_T: P_T^{{\cal A}}(\{X_i \}, Y) \rightarrow 
h_T^{*}P_T^{{\cal B}}(\{F(X_i)\}, F(Y))$ which are compatible 
with the compositions and the unit family.
It is clear how to specify this definition for the case of schemes
or analytic spaces as well as to the rational or meromorphic case.

\subsection{}

Let ${\cal C}$ be a pseudo-monoidal category over ${\cal S}$. It is called
$pseudo-braided$ if for every element $\sigma$ of the braid group $B_n$
we are given a morphism of families $\mu_{\sigma} P_T(\{X_i\},Y)\rightarrow 
P_T(\{X_{\sigma(i)}\},Y)$ which is identical on $S_T$. Here $\sigma$ acts
on $i$ as the corresponding permutation from the permutation group $S_n$.

These morphisms are required to satisfy various natural properties. Details
 can be found in [So]. Here we list them shortly:
 
 a) $\mu_{\sigma\tau}=\mu_{\sigma}\mu_{\tau}$, $\mu_1=id$ where
 $1$ denotes the unit of the group;
 
 b) for any $\sigma\in B_n$ the morphism $\mu_{\sigma}$ commutes
 with morphisms in ${\cal T}(n)$, and $\mu_1$ preserves ${\bf 1}_X$
 for any object $X$;

 c) compatibility with the composition maps (this means natural
 commutative diagram, see [So]).
 
 Now we can specify an operad ${\cal S}$ of spaces (schemes, analytic
 spaces, manifolds, etc.). Then we require that all $\mu_{\sigma}$ are
 morphisms in the corresponding category and
  arrive to various versions of pseudo-braded categories (algebraic, analytic,
  etc.). If the our spaces are irreducible (in the corresponding category)
  we can require that all $\mu_{\sigma}$ are defined in the generic point only.
  In this way we obtain the notions ot $rational$ (in the case of schemes)
  and $meromorphic$ (in the case of analytic spaces or complex manifolds)
  pseudo-braided category (see [So] for the details). The notion of a functor
  between such categories is defined in the natural way. The definitions
  in the representable case are also clear. Then we also have the notion
  of braided category, meromorphic braided category, etc. Sometimes we will
  simply call them $tensor$, $meromorphic$ $tensor$, etc.

\section{Examples}
We recall here main examples.

\vspace{2mm}

{\bf Example 1}

\vspace{2mm}

Let $G$ be a complex analytic group (the condition on $G$
can be weaker of course). We can define the following operad
of spaces. For $T \in {\cal T}(n)$ we put $S_T=G^n, G^0=id$. Then
the morphisms from the Definition 1c) are identities.
We also have an operadic composition $\gamma: G^n \times G^{k_1}\times...
\times G^{k_n} \rightarrow G^{k_1+...+k_n}$ such that
$\gamma((g_i)\times\prod_i (g_{ij}))=(g_ig_{ij})$. 

Let ${\cal C}$ be a category equipped with the action of $G$ on objects:
an object $M$ is transformed by $g \in G$ into the object called $M(g)$.

Then we say that ${\cal C}$ is a meromorphic monoidal $G$-category
if for any 
$T \in {\cal T}(n)$ we are given a functor
$\otimes_T: {\cal C}^n \rightarrow {\cal C}$ such that the 
families
$Hom_{{\cal C}}(\otimes_TX_i(g_i), Y)$ define a representable
meromorphic pseudo-monoidal structure on ${\cal C}$ such that
all morphisms from Definition 1c) are meromorphic isomorphisms.
Subsequently we have the notions of $G$-braided category or
meromorphic $G$-braided category, etc.

\vspace{2mm}

{\bf Example 2}

\vspace{2mm}

This is a specialization of the Example 1.

We take as $G$ the group ${\bf C}^{*}$ of non-zero complex numbers.
We take as ${\cal C}$ the category of finite-dimensional $U_q(g)$-
modules where $U_q(g)$ is the Drinfeld-Jimbo quantized enveloping
algebra of an affine Kac-Moody Lie algebra $g$. We call it quantum
affine algebra. The well-known fact explained in [So] is that ${\cal C}$
is a meromorphic ${\bf C}^{*}$-tensor category. We will discuss
this one and related category later in the text.

\vspace{2mm}

{\bf Example 3}

\vspace{2mm}

The usual notions of monoidal and braided (=tensor) categories are
special cases. They can be described as representable pseudo-monoidal
or pseudo-braided structures. One can take either the trivial operad
of spaces with all spaces being just one point. Or one can take
the operad of moduli of stable punctured complex curves with tangent vector 
attached to the last point. Taking a connected stratum of the real
points we get monoidal categories. In complex case we note that the
fundamental groups of the punctured curves are pure braid groups.
Thus we arrive to the description of tensor categories in terms
of local systems on punctured curves due to Deligne ([De]). See
[So] for details.

\section{Quantum affine algebras}

\subsection{}

Let $X$ be a complex manifold, $A_X$ be a bundle of Hopf algebras
on $X$ which is equipped with a flat connection $\nabla$.
The latter means that $\nabla$ has zero curvature and equalities
$\nabla(ab)=\nabla(a)b+a\nabla(b)$ and
$(\nabla\otimes 1+1\otimes \nabla)(\Delta(a))=
\Delta(\nabla(a))$ hold locally (here $\Delta$ is the comultiplication
morphism for $A_X$).

 We define a category ${\cal C}_X$ as a category of holomorphic vector
 bundles $V$ on $X$
 of finite rank, equipped with a  flat connection $\nabla_V$ which
 are $(A_X,\nabla)$-modules. This means that $V$ is a locally
 free sheaf of $A_X$-modules and the equality
 $\nabla_V(av)=\nabla_X(a)v+
 a\nabla_V(v)$ holds locally for all sections $a$ of $A_X$ and $v$ of $V$.
 Morphisms in ${\cal C}_X$ are morphisms of vector bundles compatible
 with the structures.

It is clear that  naturally defined 
kernel and cokernel of a morphism 
$(M,\nabla_M)\rightarrow (N, \nabla_N)$ belong
to ${\cal C}_X$. The tensor product is defined fiberwise, and the unit
object is defined as a trivial line bundle over $X$ equipped with
the trivial connection. This implies the following lemma.

\begin{lmm}
 
 The category ${\cal C}_X$ is an abelian monoidal category.
 \end{lmm}

Let $(M_i,\nabla_i)$ be a sequence of objects of ${\cal C}_X$,
$1\le i\le n$. We can make the tensor product $M=\boxtimes_{i=1}^{i=n}
M_i$ which is a holomorphic vector bundle on $X^n$ equipped with
the flat connection induced from tensor factors. The fiber over $(x_1,...,x_n)$
carries a structure of $\otimes_iA_{X,x_i}$-module.
Let $Z_n$ denotes infinitesimal neigbourhood of the diagonal $\{x_1=x_2=...
=x_n\}$ with the union of all diagonals $\{x_i=x_j\}$ being removed. Then we
have  functors $j_{n\ast}$ and $j_n^{\ast}$ in the category of $D$-modules
on $X^n$, where $j_n=j_{Z_n}$ is the canonical embedding of $Z_n$
into $X^n$. 

Let us consider $j_{n\ast}j_n^{\ast}(M)$. Using the flat connection on $M$
we can identify the fiber $M_{x_1,...,x_n}$ with the fiber $M_{x_1,...,x_1}$.
The latter carries the natural structure of an $A_{X,x_1}$-module. 

Let us assume that for every such a fiber and every 
permutation $\sigma\in S_n$ we are given an isomorphism of $A_{X,x_1}$-
modules $c_{\sigma}:j_{n\ast}j_n^{\ast}(\boxtimes_iM_i)\to
j_{n\ast}j_n^{\ast}(\boxtimes_iM_{\sigma(i)})$ such that:

$c_{1}=id, c_{\sigma\tau}=c_{\sigma}c_{\tau}$.

\begin{dfn} We say that we are given an infinitesimal
chiral braiding on ${\cal C}_X$ if 

a) the above-mentioned isomorphisms are functorial with respect to $M_i$;

b) they are compatible with the natural embeddings $X^n\to X^m, n\le m$
in the sense that if $\sigma \in S_m$ permutes a subsequence of $\{1,...,m\}$
consisting of $n$ elements then
$c_{\sigma}, \sigma\in S_m$ acts as $c_{\bar\sigma}\otimes id$
where $\bar\sigma$ is the corresponding element  of $S_n$.

If the above mentioned structure exists globally on the complement to
each main diagonal $\{x_1=x_2=...=x_n\}$ then we say that ${\cal C}_X$
carries a chiral braiding.

\end{dfn}

\begin{rmk} This definition can be restated in such a way that it admits 
generalization to the case of families of objects of a braided category
which are parametrized by $X$ and equipped with automorphisms of
 infinitesimally
close fibers.
Then to any planar tree with  tails numbered from $1$ to $n$, and to a
sequence of families $M_i, 1\le i\le n$ on $X$, one assigns a tensor product
$M_T=\boxtimes_TM_i$. Although it actually depends on $n$, not on $T$, 
one needs to use trees to identify the fiber $(M_T)_{x_1,...,x_n}$ with the tensor
product $\otimes_TM_{i,x_1}$ in the infinitesimal neighbourhood of the diagonal. 
Again, we get a family on $X$.
We can now define chiral braiding as before (it is also  similar
 to the definition of meromorphic braiding from 2.4 ). Notice that in this case we
 need additional 
compatibilies between trees with the same number of tails (chiral associativity)
as well as compatibilities with the gluing operation on trees.
We leave these details to the reader.

\end{rmk}

If ${\cal C}_X$ carries a chiral braided structure, then one can define meromorphic
braided structure on global sections of the objects from ${\cal C}_X$.

\subsection{}
Let $U$ be the quantized enveloping algebra of the affine 
Kac-Moody algebra $\widehat{sl(2)}$
corresponding to the parameter $q$ such that $|q|<1$.
We denote by $\overline{U}$ the quantized subalgebra $U_q(sl(2))$.
We denote by ${\cal A}$ the category of finite-dimensional $U$-modules
of type 1 (see for example [ChP]). It follows from the results of [KS] that
${\cal A}$ is a meromorphic braided category (see discussion in [So] ).

To fix the notation, we recall that $U$ is a Hopf algebra generated by
$X_i^{\pm}, K_i, K_i^{-1}, i=0,1$ subject to relations 
$$
K_iX_j^{\pm}K_i^{-1}=q^{{\pm} a_{ij}}X_j,
$$
where $(a_{ij})$ is the Cartan matrix of $\widehat{sl(2)}$,

$$
K_iK_i^{-1}=1, K_iK_j=K_jK_i,
$$

$$[X_i^{+},X_i^{-}]=(K-K^{-1})/ (q-q^{-1}),
$$
as well as quantized Serre relations for $X_i^+$ and $X_i^-$, $i=0,1$.
The coproduct $\Delta: U \to U\otimes U$ is defined by
$$\Delta(K)=K\otimes K,
\Delta(X_i^{\pm})=X_i^{\pm}\otimes K_i+K_i^{-1}\otimes X_i^{\pm}, i=0,1.
$$
Therefore the antipode is defined on generators by $S(K)=K^{-1},
S(X_i^{\pm})=-q^{{\pm}2}X_i^{\pm}, i=0,1.$

The group ${\bf C}^*$ acts on $U$ via $\phi_z(X_i^{\pm})=
z^{\pm 1}X_i^{\pm}, \phi_z(K_i)=K_i, z\in {\bf C}^*$.

\subsection{}
One of our goals will be to construct a subcategory ${\cal C}$ of ${\cal A}$
and the extend it to a category ${\cal C}_X$ of the previous subsection, 
taking $X$ to be the elliptic curve
 ${\cal E}={\bf C}^*/q^{2{\bf Z}}$.

One way to do it is to construct a category ${\cal C}$ with the following
properties:

1. ${\cal C}$ is a full monoidal rigid subcategory of ${\cal A}$;

2. ${\cal C}$ is closed with respect to taking submodules and factor
modules;

3. If $V,W \in {\cal C}$  then meromorphic brading
in ${\cal A}$ gives rise to an isomorphism $V(x)\otimes W(y)
\rightarrow W(y)\otimes V(x)$ as long as $x/ y$ does not belong
to the set $q^{2{\bf Z}}$.

Let ${\cal C}$ be  a monoidal subcategory of ${\cal A}$ which satisfies
 Property 2 only.

\begin{lmm} Property 3 holds for any two objects $V, W\in {\cal C}$
as long as it holds for simple $V$ and $W$.

\end{lmm}

$Proof$. One case use induction by $n=dimV\cdot dim W$.
The Lemma holds for $n=1$. Suppose that it holds for all $k<n$.
Let us take $V$ and $W$ such that $dimV\cdot dim W=n$.
If both objects are simple we are done. Suppose that $V$ is not
simple. Then there exists a non-trivial simple submodule $V_1
\in V$. 
 Then we have an exact sequence
$$
0\rightarrow V_1(x)\rightarrow V(x)\rightarrow V(x)/V_1(x)\rightarrow 0.
$$
We can tensor it with $W(y)$ from the left and from the right. Then
using induction assumption, functoriality
of meromorphic braiding and five-lemma we get the result. Q.E.D.

Suppose that we are given a category ${\cal C}$ which satisfies
the properties 1-3 above. Every object $V\in {\cal C}$ gives rise
to a trivial vector bundle $V_{{\bf C}^*}$ on ${\bf C}^*$ with 
the trivial connection.

We can assign to $U$ a trivial bundle $ U_{{\bf C}^*}$ of 
Hopf algebras on ${\bf C}^*$
equipped with the  connection defined by the action
 $\phi_z$ of ${\bf C}^*$.
Actually $U_{{\bf C}^*}$ is an equivariant bundle. 
Then $V_{{\bf C}^*}$ is a bundle of $ U_{{\bf C}^*}$-modules 
such that the fiber
 $ U_{{\bf C}^*,z}=U$
acts on the fiber $ V_{{\bf C}^*,z}=V$ via automorphism $\phi_z$. 

To descent these data to ${\cal E}$ we need to define automorphisms
 between
fibers at $z$ and $zq^2$ compatible with module structures.
We define them to be $id_V$ for $V_{{\bf C}^*}$ (all fibers are canonically
identified with $V$). We define an isomorphism
 $\gamma_z: U=U_{{\bf C}^*,z}
\to U_{{\bf C}^*,zq^2}=U$ as $\gamma_z(a)=\phi_{q^{-2}}(a), a\in U$.
Since $\gamma_z(\phi_z(a)v)=\phi_{zq^{-2}}(a)v=\phi_{zq^2}(\gamma_z(a))v$
for any $v\in V, a\in U$ we see that indeed all the structures are compatible and 
we obtain a bundle $U_{{\cal E}}$ of Hopf algebras on ${\cal E}$ equipped with
a connection (in fact a structure of equivariant sheaf) as well as a bundle
$V_{{\cal E}}$ of $U_{{\cal E}}$-modules.

Let $V_{{\cal E}}$ and $W_{{\cal E}}$ be two $U_{{\cal E}}$-modules.
Then $V_{{\cal E},x}\boxtimes W_{{\cal E},y}$ carries a structure
of $U_{{\cal E},x}\otimes U_{{\cal E},y}$-module. For any $x$ and $y$ 
there is an isomorphism of the fiber $U_{{\cal E},y}$ and the fiber
$U_{{\cal E},x}$ given by $\phi_{xy^{-1}}$ (to be more precise this is
the formula on ${\bf C}^*$ but it is compatible with all the structures
so it descents to the elliptic curve). This makes $V_{{\cal E},x}
\boxtimes W_{{\cal E},y}$ into $U_{{\cal E},x}\otimes U_{{\cal E},x}$--
module and via coproduct $\Delta$ into $U_{{\cal E},x}$-module.

If $(x,y)\in {\cal E}\times {\cal E} \setminus \{diag \}$ then the meromorphic
braiding in ${\cal A}$ descents to an isomorphism of $U_{{\cal E},x}$--
modules $c_{x,y}:V_{{\cal E},x}\boxtimes W_{{\cal E},y}\to
W_{{\cal E},y}\boxtimes V_{{\cal E},x}$. 

Therefore we get a chiral braiding $c_{V,W}:j_*j^*(V_{\cal E}\boxtimes
W_{\cal E})\simeq j_*j^*(W_{\cal E}\boxtimes V_{\cal E})$.

Here $j$ is the embedding of the complement of the diagonal to
${\cal E}\times {\cal E}$.

It is easy to check that in this way we have obtained a chiral braided
category (associativity constraint is trivial). Taking sections of the
bundles we get a meromorphic braided category.

\subsection{}

We are going to construct a subcategory ${\cal C}$ which satisfies  Properties 1-3 .                      
                                            
To do this we recall that for any
non-negative integer $n$ and for any non-zero complex number $a$ one has
a simple object $V_n(a)$ of ${\cal A}$. It is constructed as evaluation
representation of $U$ corresponding to the point $a$ and finite-dimensional
simple  $(n+1)$-dimensional $U_q(sl(2))$-module $V_n$. 
We call $\otimes_iV_{n_i}(a_i)$
a $standard$  module corresponding to $(a_1,a_2,...)$.
It is known (see [ChP])
that any simple object of ${\cal A}$ is a  standard one with $a_1,a_2,...$
satisfy certain properties. Namely to every $V_{n_i}(a_i)$ one assigns a finite 
set $S_{n_i}(a_i)$ called $q-string$. It is a subset of $a_iq^{2{\bf Z}}$. The condition
mentioned above says that any two $q$-strings $S_{n_i}(a_i)$ and $S_{n_j}(a_j)$
are in a $generic$ $position$ (see [ChP] for precise definitions).
If all $a_i/a_j$ do not belong to $q^{2{\bf Z}}$ then the above-mentioned strings
are in generic position.

We define ${\cal C}$ as a full rigid monoidal subcategory of ${\cal A}$ which 
is generated by $V_n(q^m), m\in 2{\bf Z}$ and closed under taking submodules
and quotients.
Clearly the trivial module ${\bf 1}$ belongs
to ${\cal C}$. Since $V_n(q^m)^*\simeq V_n(q^{2+m})$ our category is rigid monoidal.
Therefore it satisfies the Properties 1 and 2.

\begin{thm} Property 3 holds for the category ${\cal C}$.

\end{thm}
$Proof$. We will split the proof into several steps.

Step 1. Let $V=\otimes_{i+1}^{i=n}V_{k_i}(q^{l_i})$ and $W=\otimes_{j=1}^{j=m}
V_{r_j}(q^{s_j})$. If $x/y$ does not belong to $q^{2{\bf Z}}$ then two strings
$S_{k_i}(xq^{l_i})$ and $S_{r_j}(yq^{s_j})$ are in generic position
for any $i$ and $j$. Consider the tensor product $V(x)\otimes W(y)=
\otimes_iV_{k_i}(xq^{l_i})\otimes\otimes_jV_{r_j}(yq^{s_j})$. 
We can use meromorphic bradings to interchange every module from
the first group of tensor factors with every module from the second
group of tensor factors. It is easy to deduce from [ChP], sections 4 , 5 and
[KS], section 4, that meromorphic braidings do not have singularities
and hence we obtain an isomorphism of $U$-modules $V(x)\otimes W(y)
\to W(y)\otimes V(x)$.

Step 2. Let $V$ and $W$ be as on the Step 1.
If $M$ and $N$ are either both submodules or factor modules of
$V$ and $W$ then the Property 3 holds for $M$ and $N$. This is clear from the
Step 1 and the following functoriality of meromorphic braiding proven in [KS]:
if meromorphic braiding $c_{A(x),B(y)}: A(x)\otimes B(Y)\to
B(y)\otimes A(x)$ is well-defined and $f:A\to A^{\prime},
g:B\to B^{\prime}$ are morphisms of $U$-modules then
$(f\otimes g)(c_{A(x),B(y)})$ is a well-defined isomorphism 
$A^{\prime}(x)\otimes B^{\prime}(y)\to B^{\prime}(y)\otimes A^{\prime}(x)$.

Similarly one can prove the Property 3 in case if  $M$ is a submodule of $V$
and $N$ is a factor module of $W$. 
Using functoriality of meromorphic braiding once again, we can prove the Property 3
for any object of ${\cal C}$.
Q.E.D.

\subsection{}
It is natural to ask for a description of simple objects of ${\cal C}$.
 We start with the following elementary  result.

\begin{thm} Any simple submodule or factor module of $V=\otimes_{i=1}^{i=n}
V_{k_i}(q^{l_i})$ is of the form $M=\otimes_{i=1}^{i=n}V_{r_i}(q^{s_i}).$

\end{thm}
$Proof$. Let us use induction by $k=\sum_ik_i$. For $k=0$ the result is
obvious. Suppose it holds for all $k<m$. Let us prove it for $k=m$.
We prove it for submodules
only. The case of factor modules easily follows if we take duals.

If $V$ is simple we have nothing to prove.
Let $M\subset V$ be a non-trivial submodule. Since $V$ is not simple
there are two $q$-strings $S_{k_i}(q^{l_i})$ and $S_{k_j}(q^{l_j})$ which
are not in generic position. We may assume that $i=1, j=2$ (otherwise
use the bradings to get two consequitive strings  not in generic
position). Then according to  [ChP], Section 4.9, there is a unique  simple
submodule of $X=V_{k_1}(q^{l_1})\otimes V_{k_2}(q^{l_2})$ of the type
$A=V_d(q^e)\otimes V_f(q^h)$ with the factor module $B$ of similar type.
Then we have an exact sequence
$$
0\to A\otimes_{i \ge 3}V_{k_i}(q^{l_i})\to V\to B\otimes_{i\ge 3}V_{k_i}(q^{l_i})
\to 0$$

If simple module $M$ intersects $A\otimes_{i\ge 3}V_{k_i}(q^{l_i})$
then  it belongs to it and
the result follows from the induction assumption. 
If $M$ does not intersect this submodule
then it is projected isomorphically to $B\otimes_{i\ge 3}V_{k_i}(q^{l_i})$
and again the result follows by induction. Q.E.D.

\vspace{2mm}

On ther other hand one can try to use the notion of $q$-character
introduced in [FR] in order to describe the Grothendieck ring of
the category ${\cal C}$.
For example it is natural to expect an affirmative answer to the following

{\bf Question} Let $\chi_q$ be the $q$-character (notation from [FR]).
Is it true that the Grothendieck ring $K_0$ of our category ${\cal C}$ is isomorphic
to the subring of the Grothendieck ring of $U$ generated by $t_{q^n}$,
where $t_{q^n}$ is the class of $V_1(q^n), n\in 2{\bf Z}$?

\vspace{2mm}

Clearly $K_0({\cal C})$ contain the subring ${\bf Z}[t_{q^n}], n\in 2{\bf Z}$.
Thus the question is whether the $q$-character of an object from ${\cal C}$
belongs to this ring.

As was pointed to me by Ed Frenkel (private communication) the answer to the Question 
is positive, and it can be generalized
to the higher rank case.
We sketch his arguments below in the case of quantum affine algebra
of $sl(n)$. We denote by $\Gamma_n$ the set $2{\bf Z}$ if $n=2$ and
the set ${\bf Z}$ if $n>2$.

Let us consider the monoidal category ${\cal C}_n$ generated by evaluation 
representations
 $V_{\omega_i}(q^l)$ of $U_q(\widehat {sl(n)})$ where $l\in \Gamma_n$ and $\omega_i$
is the fundamental weight of $U_q(sl(n))$ $1\le i\le n-1$. 
Then the $q$-character of the tensor product of such
representations belongs to $A_n={\bf Z}[t_{i,q^l}]$ where $t_{i,q^l}$ is the class of 
$V_{\omega_i}(q^l)$ in the representation ring of $U_q(\widehat {sl(n)})$, and
$l \in \Gamma_n$. Thus for $n=2$ we have ${\cal C}_2={\cal C}$
and $A:=A_2$ is  expected to be isomorphic to  $K_0({\cal C})$.

\begin{thm} A simple object in ${\cal C}_n$ is isomorphic to a subquotient of
a tensor product $\otimes_iV_{\omega_i}(q^{l_i}), l_i\in \Gamma_n$. 

In the case $n=2$ it
is isomorphic to  a tensor product $\otimes_iV_{n_i}(q^{l_i})$
 where $V_{n_i}(q^{l_i})$ is a standard module, $l_i \in {\bf Z}$.

\end{thm}

$Proof$. According to Chari and Pressley a simple $U_q(\widehat{sl(n)})$-module
is isomorphic to a subquotient of $\otimes_iV_{\omega_i}(a_i)$ (we have an isomorphism
to  a tensor product of standard modules for $n=2$).
 Then one uses the fact that the $q$-character of
the simple module contains the dominant term (terminology and notation from [FR],
Section 4)
equal to $\prod_iY_{i,a_i}$, and is the sum with positive coefficients of monomilas
in $Y_{i,a_iq^{n_i}}^{\pm 1}$. This implies that all $a_i\in q^{{\Gamma_n}}$. Q.E.D.

Let $V=\otimes_iV_{\omega_i}(q^{l_i})$.

\begin{thm} The $q$-character of every subquotient of $V$ belongs
 to $A_n$. 
\end{thm}
$Proof$. Let $L_j, 1\le j\le m$ be the set of simple objects which appears in
the composition series of $V$. Then $\chi_q(V)=\sum_j \chi_q(L_j)$.

Every $L_j$ is  the highest weight module over $U_q(\widehat{sl(n)})$.
The $q$-character of such a module was computed in [FR]. 
It is equal to
a sum of monomials in $t_{i,a}$ with positive coefficients.
It follows from the previous theorem that if $t_{i,a}$ appears in such monomials we have
$a\in q^{{\Gamma_n}}$. Hence $\chi_q(L_j)\in A_n$ for every $j$. 
This implies the theorem. Q.E.D.

\begin{cor} The Grothendieck ring $K_0({\cal C}_n)$ is isomorphic to $A_n$.

\end{cor}

$Proof$. It is easy to see that $K_0({\cal C}_n)$ is generated by the isomorphism
classes of subquotients of tensor products of the type
$V=\otimes_iV_{\omega_i}(q^{l_i})$ and then apply the previous theorem. Q.E.D.

There is little doubt about positive answer to the Question for an arbitary quantum
affine algebra. One can try to use this fact 
in order to construct  the corresponding 
meromorphic braided category on an elliptic curve.

\vspace{30mm}

{\bf References}

\vspace{3mm}

[BD] A. Beilinson, V.Drinfeld, Chiral algebras, preprint, 1995.
\vspace{2mm}

[B] R. Borcherds, Vertex algebras, q-alg 9706008.

\vspace{2mm}

[ChP] V. Chari, A. Pressley, Quantum affine algebras, Comm. Math. Phys.,
142(1991), 261-283.

\vspace{2mm}

[FR] E. Frenkel, N. Reshetikhin, The $q$-characters of representations of quantum
affine algebras and deformations of $W$-algebras, q-alg 9810055.

\vspace{2mm}

[GK] V. Ginzburg, M. Kapranov, Koszul duality for operads, Duke math. J.,76:1(1994),
203-272.

\vspace{2mm}

[KM] M. Kontsevich, Yu. Manin, Gromov-Witten classes, quantum cohomology and
enumerative geometry, Comm. Math. Phys., 164(1994),525-567.

\vspace{2mm}

[KS] D. Kazhdan, Y. Soibelman, Representations of quantum affine algebras,
Selecta Math., New Series, 1:3(1995), 537-595.

\vspace{2mm}

[L] T. Leinster, General operads and multicategories, CT 9810053.

\vspace{2mm}

[La] J. Lambeck, Deductive systems and categories, Lect. Notes Math.,
86(1969), 76-122.

\vspace{2mm}

[So] Y. Soibelman, Meromorphic tensor categories, q-alg 9709030.

\end{document}